\documentclass[12pt]{amsart}
\usepackage{amsthm,amstext,amsmath,amscd,amssymb,latexsym,mathrsfs}

\usepackage{color}
\usepackage{amsfonts,array}
\usepackage{todonotes}
\usepackage{verbatim}
\usepackage{float}	
\usepackage{pgf,tikz}
\usepackage[toc,page]{appendix}
\usepackage[all,cmtip]{xy}
\usepackage{stmaryrd}
\usepackage{xparse}
\usepackage{listings}
\usepackage{enumitem}
\setlist{nosep, left=0pt}

\usepackage{tabularx,booktabs}

\usepackage{ragged2e} 
\usepackage{needspace}   
\usepackage{ltablex}   
\keepXColumns         
\usepackage{placeins} 

\newcolumntype{L}{>{\RaggedRight\arraybackslash}p{3.5cm}}
\newcolumntype{Y}{>{\RaggedRight\arraybackslash}X}

\usepackage{caption}
\usepackage{subcaption}
\usepackage[T1]{fontenc}
\usepackage{lmodern}
\usepackage{etoolbox}
\usepackage{xcolor}

\usepackage[shortalphabetic]{amsrefs}
\usepackage[all]{xy}
\usepackage{newtxtext}
\usepackage{newtxmath}
\usepackage{framed}
\usepackage{graphicx}
\usepackage{wasysym}
\usepackage{pdfpages}

\usepackage[thinlines,thiklines]{easybmat}

\usepackage{microtype} 

\usepackage{mathtools}  

\usepackage{bbm}
\usepackage[bookmarks=true, colorlinks=true, linkcolor=blue!50!black,
citecolor=GoetheBlue, urlcolor=orange!50!black, pdfencoding=unicode]{hyperref}
\usepackage[left=3.5cm,right=3.5cm,top=3cm,bottom=3cm,footskip=1cm
]{geometry}
\makeatletter
\def\thm@space@setup{%
  \thm@preskip=5pt \thm@postskip=5pt
}
\makeatother

\makeatletter
\renewcommand\@biblabel[1]{[#1]\quad}%

\@ifpackageloaded{amsrefs}{%
  \ifcsname BibLabel\endcsname
    \usepackage{etoolbox}
    \patchcmd{\BibLabel}{\hfill}{}{}{}%
    \renewcommand{\BibLabel}{%
      \ifcsname Hy@raisedlink\endcsname
        \Hy@raisedlink{\hyper@anchorstart{cite.\CurrentBib}\hyper@anchorend}%
        [\thebib]\quad
      \else
        [\thebib]\quad
      \fi
    }%
  \fi
}{}%
\makeatother

\def\Bl{\operatorname{Bl}}

\numberwithin{equation}{section}

\newcommand{\equ}{\ensuremath{\,=\,}}

\newcommand{\deq}{\ensuremath{\stackrel{\textrm{def}}{=}}}

\DeclareMathOperator{\lint}{int}

\newcommand{\Ra}{`\ensuremath{\Rightarrow}'  }
\newcommand{\La}{`\ensuremath{\Leftarrow}'  }
\newcommand{\LI}{`\ensuremath{\subseteq}'  }
\newcommand{\RI}{`\ensuremath{\supseteq}'  }

\newcommand{\BA}{{\mathbb{A}}}

\newcommand{\BC}{{\mathbb{C}}}

\newcommand{\BG}{{\mathbb{G}}}

\newcommand{\BN}{{\mathbb{N}}}

\newcommand{\BP}{{\mathbb{P}}}
\newcommand{\BQ}{{\mathbb{Q}}}
\newcommand{\BR}{{\mathbb{R}}}

\newcommand{\BZ}{{\mathbb{Z}}}

\newcommand{\CC}{{\mathcal C}}

\newcommand{\CL}{{\mathcal L}}

\newcommand{\CW}{{\mathcal W}}

\DeclareMathOperator{\Spec}{Spec}
\DeclareMathOperator{\Proj}{Proj}

\DeclareMathOperator{\rank}{rank}

\DeclareMathOperator{\Pic}{Pic}
\DeclareMathOperator{\Cl}{Cl}

\DeclareMathOperator{\Eff}{Eff}
\DeclareMathOperator{\Nef}{Nef}

\DeclareMathOperator{\Cone}{Cone}
\DeclareMathOperator{\Mov}{Mov}

\DeclareMathOperator{\NE}{NE}

\DeclareMathOperator{\Cox}{Cox}

\DeclareMathOperator{\Rel}{Rel}
\DeclareMathOperator{\Hom}{Hom}

\newcommand{\Stab}{\mathop{\rm Stab}\nolimits}

\newcommand{\Gen}{{\rm Gen}}

\definecolor{GoetheBlue}{RGB}{0,97,143}
\usepackage[bookmarks=true, colorlinks=true, linkcolor=blue!50!black,
citecolor=GoetheBlue, urlcolor=orange!50!black, pdfencoding=unicode]{hyperref}
\usepackage[left=3.5cm,right=3.5cm,top=3cm,bottom=3cm,footskip=1cm
]{geometry}
\makeatletter
\def\thm@space@setup{%
  \thm@preskip=5pt \thm@postskip=5pt
}
\makeatother

\newtheorem*{theorem*}{Theorem}

\newtheorem{theorem}{Theorem}[section]
\newtheorem{lemma}[theorem]{Lemma}

\newtheorem{proposition}[theorem]{Proposition}
\newtheorem{corollary}[theorem]{Corollary} 

\theoremstyle{definition}
\newtheorem{definition}[theorem]{Definition}
\newtheorem{example}[theorem]{Example}

\newtheorem{remark}[theorem]{Remark}

\theoremstyle{plain}
\newtheorem{introthm}{Theorem}

\def\phi{\varphi}
\def\epsilon{\varepsilon}
\def\setminus{\smallsetminus}
\let\oldbullet\bullet
\def\bullet{{\mathchoice{\oldbullet}%
                        {\oldbullet}%
                        {\scriptscriptstyle\oldbullet}%
                        {\oldbullet}}}
\let\oldemptyset\emptyset
\let\emptyset\varnothing

\title{Stability Conditions for Multigraded Rings} 

\author{Felix G\"obler}

\begin{document}






\maketitle

\begin{abstract}
Let $D$ be a finitely generated abelian group and $S$ a $D$-graded ring. We introduce a notion of geometric (semi)stability for the action of the diagonalizable $S_0$-group scheme $G = \Spec(S_0[D])$ on the $S_0$-scheme $X = \Spec(S)$, and show that the multihomogeneous spectrum $\Proj^D(S)$ is the geometric quotient of the open subset of geometrically semistable points $X^{\mathrm{gss}}
\subseteq X$ by $G$.

Applying this construction to the Cox ring of a toric (pre)variety (resp. of a Mori dream space) yields new insights both into the quotient construction of such varieties and into the computation of those quotients, as well as into the secondary fan (resp. into the Mori chamber decomposition).
In particular, we give an equivalent formulation for Mori equivalence in terms of certain homogeneous sections of the Cox ring.
\end{abstract}

\tableofcontents


\section*{Introduction}
Algebraic geometry deals with the interconnections of algebra and geometry. One of the first instances is the anti-equivalence of the category of affine schemes and the category of rings. But affine geometry has its limitations, most notably the absence of completions. 

Projective varieties, on the other hand, are proper, but the ring of global sections is constant. However, there still is a coordinate ring $S$, graded by the natural numbers, such that the basic affine opens of this projective variety are given in terms of degree-zero localizations $\Spec(S_{(f)})$ for homogeneous $f \in S$.

In general, coordinate rings of varieties might carry gradings by finitely generated abelian groups.
In geometry, for example, such rings appear very naturally as Cox rings of toric (pre-)varieties $X$ (and more generally of Mori dream spaces), where the grading is given by the divisor class group $\Cl(X)$.

For toric varieties $X_\Sigma$, Cox showed that the Cox ring together with the ideal defined by the maximal cones of the fan $\Sigma$ recovers the toric variety as a GIT quotient (see \cite{Cox}, Theorem 2.1).
Applying Brenner--Schröers' notion of relevance to Cox rings yields an intrinsically defined algebraic invariant, namely the \emph{irrelevant ideal} $S_+$, which is the ideal generated by all relevant elements. The question is, what information does it contain? Can we describe this information geometrically?

Concretely, let $D$ be a finitely generated abelian group.
Given a $D$-graded ring $S$, $X= \Spec(S)$, the set of \emph{semistable} points associated to $d \in D_\BR$ is defined as
        \begin{align}\label{eq:toric_ss}
            X^{\text{ss}}(d) \deq \{ x \in \Spec(S) \mid \exists n > 0, f \in S_{nd}: f(x) \neq 0\} \subseteq \Spec(S).
        \end{align}

However, the notion of semistability does not suffice to ensure that the quotient $X^{\text{ss}}(d)\sslash G$ is of the correct dimension $\dim(\Spec(S)) - \rank(D)$ (cf.\ \cite{paper1}, Lemma 2.20).

We show that the right notion of semistability for the action of $G $ on $\Spec(S)$ has a geometric component. 
Specifically, a point $x \in X$ must have the property that
$$\sigma(x) := \overline{\Cone}(d \in D \mid \exists f \in S_d: f(x) \neq 0)$$ has maximal dimension in $D_\BR$.
Equivalently, the $G$-orbit of $x$ has dimension equal to the rank of $D$.
This way, we can ensure that the quotient has the desired dimension. In fact, it forces the homogeneous $f \in S$ such that $f(x)\neq 0$ to be relevant. We call this property \emph{geometric semistability} and prove the following new result that applies to all multigraded rings.

\begin{introthm}[\autoref{thm:X^gss_as_relevant_cone}]
    Let $D$ be a finitely generated abelian group, $S$ a $D$-graded ring and $G = \Spec(S_0[D])$. 
     The set of geometrically semistable points with respect to the action of the diagonalizable 
$S_0$-group scheme $G$ on $X=\Spec(S)$ is given by
    \begin{align*}
        X^\text{gss} \equ \bigcup_{f \in \Rel^D(S)} D(f) \equ D(S_+).
    \end{align*}
    In particular, $\Proj^D(S)$ is characterized in terms of $\Proj^D(S) := X^\text{gss} / G$, so it can be described without the notion of relevance. 
\end{introthm}

 Note that by definition, every stable point is geometrically semistable. Hence, stable subsets for the action of $G$ on $\Spec(S)$ correspond to maximal separated subsets of $X^\text{gss}$, i.e.\ there may be more than one choice.

From now on in this introduction, we assume that $S_0$ is an algebraically closed field $k$ of characteristic zero.
If we specialize to the case when $S$ is the Cox ring of a quasiprojective simplicial toric variety $X_\Sigma$ (over $k$), these choices of subsets of $X^\text{gss}$ are related to the \emph{VGIT} (variation of geometric invariant theory) of $X_\Sigma$, that is, the study of how GIT quotients $\Spec(\Cox(X_\Sigma)) //_\chi \Spec(\Cox(X_\Sigma)_0[D])$ change as the linearization/character $\chi$ varies (see \cite{T} for more details). In particular, the stable subsets correspond to linearizations of birational models of $X_\Sigma$. 
Thus, using the data of a pair $(S, B)$, where $S$ is a multigraded polynomial ring\footnote{We use the terms \emph{$D$-graded} and \emph{multigraded} interchangeably. The latter term is used when the explicit grading group is clear from the context.} with finitely many variables and $B$ is an ideal in $S$ generated by a subset of the relevant elements, we can describe the stable subsets of $\Proj^D(S)$ with simplicial toric varieties.

From this point of view, it is a very natural question to ask how the $D$-graded Proj construction of the Cox ring $S$ of a simplicial normal toric (pre-)variety $X_\Sigma$ is related to $X_\Sigma$. 
Let $X_\Sigma$ be a quasiprojective simplicial toric variety\footnote{\textcolor{black}{We work with normal toric varieties over an algebraically closed field $k$ of characteristic zero.}} such that $\Sigma$ has full-dimensional convex support (i.e. $X_\Sigma$ is projective over an affine toric variety of the same dimension).
In this case $X$ has a \emph{secondary fan} by \cite{CLS}, §14.
In the toric setting, the VGIT of $X_\Sigma$ is known to be given in terms of a chamber decomposition of the weight space.

Specifically, the variation of the birational models of $X$ is described in terms of the chambers of its secondary fan. Each chamber corresponds to a generic character $\chi = \chi^a \colon G \to \BG_m$, i.e.\ a character $\chi^a$ such that $a$ lies in the interior of a maximal-dimensional cone of the secondary fan. In turn, every generic character defines an irrelevant ideal $B(\chi^a)$ (see Defintion~\ref{def:toric_irrel_ideal}) in $S = \Cox(X_\Sigma)$, where $a \in \Cl(X_\Sigma)$  and $\Cl(X_\Sigma)$ is a finitely generated abelian group (cf.\  \cite{CLS}, §14). 
Note that this characterization does not apply solely to toric varieties (see \cite{KKL}, Theorem 5.2).

A direct benefit of the Brenner--Schröer Proj construction is that we can characterize generic characters using relevant elements. This is due to the fact that the Cox ring of a toric variety is always a $\Cl(X)$- graded polynomial ring with finitely many variables over $k$.

\begin{introthm}[\autoref{thm:secondary_relevant}]
Let $X$ be a quasiprojective simplicial toric variety (over $k$) such that $\Sigma$ has full-dimensional convex support, $S = \Cox(X)$, $\Gen^D(S)$ is a generating system of $S_+$, and let $\chi = \chi^a \in \widehat{G}$ be a character. Then $\chi^a$ is generic if and only if
    \begin{align*}
        B(\chi^a) \equ ( f\in \Gen^D(S) \mid a \in \CC_D(f) ) \unlhd S.
    \end{align*}
    In particular, the chamber $\sigma_a$ in the secondary fan containing $a$ is given by
    \begin{align*}
           \sigma_a \equ \bigcap_{\substack{f \in \Gen^D(S) \\ a \in \CC_D(f)}} \CC_D(f).
    \end{align*}
\end{introthm}

Also note that, conversely, any multigraded polynomial ring with finitely many variables over $k$ yields such a chamber structure by \cite{paper4}, Theorem A.

While the proofs rely on translating known arguments into the multigraded framework of Brenner–Schröer, the resulting perspective is genuinely new. It not only provides motivation for investigating the $D$-graded Proj construction, but also yields a clearer structural understanding of the construction.
A new result is the description of models of the secondary fan through the $D$-graded Proj associated to its chamber:

\begin{introthm}[\autoref{thm:list_equiv_chamber_var}]
        Let $X$ be a quasiprojective simplicial toric variety (over $k$) such that $\Sigma$ has full-dimensional convex support, $S = \Cox(X)$ is $\Cl(X)$-graded, and let $\chi = \chi^a$ be a generic character. Then the following objects coincide:
    \begin{enumerate}[label=(\roman*)]
        \item $X_\Sigma$, where $\Sigma$ is the normal fan of $P_a$.
        \item $\Proj^\BN(R_{\chi^a})$.
        \item $\Proj^{\Cl(X)}(\sigma_a \cap \Cl(X)) = \Proj^{\Cl(X)}(\bigoplus_{d \in \sigma_a \cap \Cl(X)} H^0(X, d))$.
        \item $\Proj^{\Cl(X)}(\BC[h_1, \ldots, h_r])$, where $\deg(h_i)$ are the ray generators of $\sigma_a$.
        \item $\Proj^{\Cl(X)}_B(S)$, for $B = B(\chi^a)$.        
    \end{enumerate}
\end{introthm}
We also prove the corresponding statements for Mori dream spaces (see Theorem~\ref{thm:MDS_relevant} and Corollary~\ref{cor:MDS_chamber=sectionring}).

This paper is organized as follows: In Chapter 1, we recall some basic definitions and properties on multigraded rings and the corresponding Proj construction from \cite{BS} and our previous paper \cite{paper1}, and give a separation criterion that motivates Chapter 3. In §2, we describe $\Proj^D(S)$ in terms of stability conditions, which is independent of the notion of relevance. In Chapter 3, we apply our insights from §2 to Cox rings of toric (pre)varieties, and prove that the stable subsets of $\Proj^D(S)$ are given in terms of the chambers of its secondary fan.
In  Chapter 4, we apply the results from §2 to Mori dream spaces.


{\normalfont\bfseries Acknowledgements.}
First, I would like to thank my supervisor, Alex Küronya, for his tremendous support and insightful feedback.
I am also very grateful to Johannes Horn, Andrés Jaramillo Puentes, Kevin Kühn, Jakob Stix, Martin Ulirsch, and Stefano Urbinati for many helpful discussions and valuable comments.
In particular, I want to thank an anonymous referee for his suggestions.

Partially funded by the Deutsche Forschungsgemeinschaft (DFG, German Research Foundation) TRR 326 \textit{Geometry and Arithmetic of Uniformized Structures}, project number 444845124, and by the LOEWE grant \emph{Uniformized Structures in Algebra and Geometry}.



\section{Relevant Cones}\label{sec:rel_cones}
Let $D$ be a finitely generated abelian group and $S$ a multigraded ring. A homogeneous $f \in S$ is called \emph{relevant} if the group of homogeneous units $D^f$ of $S_f$ has finite index in $D$.
The \emph{irrelevant ideal} $S_+$ is defined to be the ideal generated by all relevant elements in $S$. 
By \cite{BS}, Lemma 2.4, $S$ is noetherian if and only if $S_0$ is noetherian and $S$ is an $S_0$-algebra of finite type.
We denote by \(\Gen^D(S)\) the finite, essentially canonical generating set of \(S_+\) provided by \cite{paper1}, Lemma 1.19 (if $S$ is noetherian).
For basic properties of multigraded rings and several equivalent notions for relevance, we refer to our previous paper \cite{paper 1}, as well as to \cite{BS}, §2, \cite{KU}, §2, and \cite{May}, §2.2.
In particular, the relevance of $f$ is equivalent to the \emph{weight cone}
\begin{align*}
    \CC_D(f) &\deq \overline{\Cone}( \{\deg_D(g) \mid g \text{ is a homogeneous divisor of $f^k$ for some $k\ge 0$}\}) 
\end{align*}
of $f \in S_d$ having maximal-dimension in $D_\BR$. 
We want to address some properties related to the relevant (weight) cones. Then we will prove that relevant cones determine the chamber structure of the weight space 
\begin{align*}
    \sigma(S) \equ\overline{\Cone}(d \in D \mid S_d \neq 0) \subseteq D_\BR.
\end{align*}
Note that we overall assume the $D$-grading on $S$ to be \emph{effective}, so that $D$ is generated by its support. (cf.\ \cite{paper1}, Definition 1.10).

First, we want to give a convex geometric criterion for separatedness of (open subsets of) $\Proj^D(S)$. It is a stronger version of \cite{BS}, Proposition 3.3, as we also show $(a) \implies (b)$, but only for factorially graded rings (where homogeneous elements have unique homogeneous factorization).
Note that separatedness is equivalent to the surjectivity of all multiplication maps $\mu_{(fg)}\colon S_{(f)} \otimes S_{(g)} \to S_{(fg)}$ for all relevant $f, g \in S$ (cf.\ \cite{BS}, proof of Proposition 3.3).

\begin{proposition}\label{prop:open_sep}
Let $S$ be a factorially graded integral domain such that $\deg(S^\times) = \{0\}$, and take $f\neq g \in \Gen^D(S)$. Then the following are equivalent: 
    \begin{enumerate}[label=(\alph*)]
        \item $\mu_{(f g)}$ is surjective.
        \item $\lint\left(\CC_D(f) \cap \CC_D(g)\right) \ \neq \ \emptyset$.
    \end{enumerate}
\end{proposition}

\begin{proof}
    Let $\sigma_1 := \CC_D(f)$ and $\sigma_2 = \CC_D(g)$.
By \cite{BS}, Proposition 3.3 we already know that $(b) \implies (a)$. 
\vspace{0.5em}

$(a) \implies (b)$: Assume that $\lint(\sigma_1 \cap \sigma_2) = \emptyset$. In particular, we know that $d=\deg(f) \not\in \sigma_2$ and $e = \deg(g) \not\in \sigma_1$ (since otherwise we may assume $d \in \partial \sigma_1 \cap \lint(\sigma_2)$), yielding $\lint(\sigma_1 \cap \sigma_2) \neq \emptyset$. Thus, let $d \not\in\sigma_2$. Then there exists at least one homogeneous divisor $h$ such that $h \mid g^k$ for some $k \in \BN$ and $h \not\mid f^l$ for all $l$. For the same reason, there is a homogeneous divisor $h'$ such that $h' \mid f^l$ for some $l \in \BN$ and $h' \not\mid g^m$ for all $m$. 
Note that the elements $h, h'$ exist by assumption. As $S$ is factorially graded with units in $S_0$, the non-existence of $h$ and $h'$ would imply that $f = g$.
Now $hh'$ is invertible in $S_{(fg)}$, but not in $S_{(f)}$ or $S_{(g)}$. Thus $\mu_{(f g)}$ cannot be surjective. 
\end{proof}

In particular, we get separated subsets of $\Proj^D(S)$ for each choice of an ideal $B \unlhd S$ generated by some subset of the relevant elements, such that the intersection of all $\CC_D(f)$ for $f \in B$ is maximal-dimensional. If $S = R[T_1,\ldots, T_n]$ for some ring $R$, this result can be visualized geometrically, using the exact sequence  
\begin{align}
    0 \ \to \ M \ \to \ \BZ^n \ \stackrel{\gamma}{\to} \ D \ \to 0,
\end{align}
where $\gamma$ sends $e_i$ to $\deg(T_i) \in D$ and $M$ is its kernel.
Since we can dualize this exact sequence using $\Hom(-,\BZ)$, the relevant elements define cones $\sigma_f$, dual to $\CC_D(f)$, and each maximal choice of a separated subset corresponds to a simplicial fan, constructed by the $\sigma_f$ (cf.\ \cite{BS}, Proposition 3.4).
If a union of cones $\sigma_f$ forms a fan, the intersection of the dual cones $\CC_D(f)$ is maximal-dimensional.
This duality is called \emph{Gale duality} (cf.\ \cite{CLS}, §14). We will be more precise about the details in Chapter~\ref{sec:VGIT}.

The following Corollary is equivalent to \cite{KU}, Remark 2.11.

\begin{corollary}
Let $\rank(D) \ge 2$.
Then the degrees of relevant elements $f \in S$ cannot lie on the boundary of $\sigma(S)$, where $\sigma(S)$ is given by the union of all $\CC_D(f)$ for homogeneous $f \in S$.
\end{corollary}

\begin{proof}
As $\deg(f) \in \lint(\CC_D(f))$ and $\CC_D(f) \subseteq \sigma_S$, $\deg(f)$ has to be in the interior of $\sigma_S$.
\end{proof}

\begin{remark}
    Let $X_\Sigma$ be a quasiprojective simplicial toric variety (over $k$) such that $\Sigma$ has full-dimensional convex support. Then the above Corollary states that the degrees of relevant elements from $\Cox(X_\Sigma)$ lie in the interior of the secondary fan of $X_\Sigma$, which lives in $\Cl(X_\Sigma)_\BR$ with support $\Eff(X_\Sigma)$.
\end{remark}

We can also show that the relevant cones are rational and polyhedral when $S$ is noetherian and factorially graded with units in $S_0$.

\begin{lemma}\label{lem:polyhedral}
Let $S$ be a noetherian factorially $D$-graded ring with $\deg(S^\times) = \{0\}$ and $f \in S$ relevant. Then $\CC_D(f)$ is a rational polyhedral cone.
\end{lemma}

\begin{proof}
    First, let $d \in C^f$. Then there is some $m \in \BN$ and $g \in S_d$ such that $g \mid f^m$. By \cite{BS}, Lemma 2.4, $S$ is an $S_0$-algebra of finite type. Since $f$ is relevant, there exist homogeneous $h_j$ and non-negative integers $\epsilon_j$ such that $f = h_{i_1}^{\epsilon_1} \cdot \ldots \cdot h_{i_r}^{\epsilon_r}$ (we ignore the $S_0$-factors from possible units). Thus there are $\epsilon_j' \in \BZ_{\ge 0}$ such that $g = h_{i_1}^{\epsilon_1'} \cdot \ldots \cdot h_{i_r}^{\epsilon_r'}$. We deduce that $d = \sum_{j=1}^r \epsilon_j' \cdot \deg(h_{i_j})$. Hence $\CC_D(f) = \overline{\text{Cone}}(C^f) = \overline{\text{Cone}}(\{\deg(h_{i_1}), \ldots, \deg(h_{i_r})\})$ is polyhedral. Since $\deg(h_{i_j}) \in D$, $\CC_D(f)$ is even rational.
\end{proof}

\begin{remark}\label{rem:matroid}
In the setting of Lemma~\ref{lem:polyhedral}, let's think of $D$ as a configuration of vectors in $D_\BR$. We see that $\Gen^D(S)$ corresponds to the maximal linearly independent subsets of that configuration, i.e.\ it corresponds to the set of bases of a matroid.
\end{remark}


\section{Stability}
There are different approaches to semistability and stability in the literature that can be applied to our situation.
Both are essentially induced by \cite{MFK}, Definition 1.7.
On the one hand, \cite{CRB} define a notion of semistability for a quasitorus acting on an affine scheme in Definition 3.1.2.1. 
While \cite{CRB} work with graded Cox rings of certain types of schemes, we transport the definition to arbitrary multigraded rings.

\begin{definition}
Let $S$ be a $D$-graded ring, $X = \Spec(S)$, and $G = \Spec(S_0[D])$. The set of \emph{semistable} points associated with an element $d \in D_\BQ$ is defined as
    \begin{align*}
        X^\text{ss}(d) \deq \{x \in X \mid \exists n>0, f \in S_{nd}: f(x) \neq 0\}.
    \end{align*}
\end{definition}

On the other hand, in the context of toric varieties, \cite{CLS} define this set using characters (i.e.,\ linearized line bundles).

\begin{definition}[\cite{CLS}, Definition 14.1.6]
Let $S$ be a $D$-graded polynomial ring with finitely many variables over an algebraically closed field of characteristic zero.
    Let $\chi\colon G\to\BG_m$ be a character, $X = \Spec(S)$, and $G = \Spec(S_0[D])$. A point $x \in X$ is called \emph{semistable}, if there are $d > 0$ and $s \in \Gamma(X, \CL_{\chi^d})^G$ such that $s(x) \neq 0$. The point is called \emph{stable} if additionally the isotropy subgroup $G_x$ is finite, and all $G$-orbits in $D(s)$ are closed in $D(s)$.
\end{definition}

Both notions of semistability are applicable and even coincide for polynomial rings $S$ with finitely many variables, since 
\begin{align*}
    S_0[D] \equ \bigoplus_{d\in D} S_0 \chi^d ,
\end{align*}
\begin{align}\label{eq:1}
    f \in S_d \ \iff \ f(gx) \equ \chi^d(g) f(x) \ \ \text{ for all } g \in G, x \in  \Spec(S),
\end{align}
and therefore \(\Gamma(X,\mathcal L_{\chi^d})^G \cong S_d\).
We refer to \cite{paper1}, Section 1.3 and Appendix A, for more details on this.

Although we generally work with a multigraded ring $S$, the following notion of semistability applies to the Cox ring of any normal integral scheme $X$ whose divisor class group $\Cl(X)$ is finitely generated.

The goal of this section is to identify the set of geometrically semistable points with the relevant locus $D(S_+)$.

\begin{definition}\label{def:semi_stable}
Let $S$ be a $D$-graded ring and $x \in X = \Spec(S)$.
    \begin{enumerate}[label=(\arabic*)]
        \item We define the \emph{orbit cone} $\sigma(x) \subseteq \sigma(S)$ of $x \in \Spec(S)$ to be 
        \begin{align*}
            \sigma(x) \deq \overline{\Cone}(d \in D \mid \exists f \in S_d: f(x) \neq 0) \subseteq D_\BR.
        \end{align*}

        \item Following \cite{CRB}, Definition 3.1.2.6, the \emph{GIT cone} of an element $d \in \sigma(S)$ is defined as the (nonempty) intersection of all 
        orbit cones containing it:
        \begin{align*}
            \lambda(d) \deq \bigcap_{x\in\Spec(S), d \in \sigma(x)} \sigma(x).
        \end{align*}

        \item The set of \emph{semistable} points associated to $d \in D_\BQ$ is defined as
        \begin{align*}
            X^{\text{ss}}(d) \deq \{ x \in \Spec(S) \mid \exists n > 0, f \in S_{nd}: f(x) \neq 0\} \subseteq \Spec(S).
        \end{align*}

        \item The set of \emph{geometrically semistable}  points associated to $d \in D_\BQ$ is defined as
        \begin{align*}
            X^{\text{gss}}(d) \deq   \{x \in \Spec(S) \mid x \in X^\text{ss}(d) \text{ and } \dim(\sigma(x)) = \rank(D)\}.
        \end{align*}

        \item The set of geometrically semistable points associated to the action of $G$ on $\Spec(S)$ is defined to be $X^{\text{gss}} = \cup_{d \in D} X^{\text{gss}}(d)$.

        \item We say that $x$ is a \emph{stable}\footnote{By \cite{paper1} Theorem A, the quotient morphism $\Spec(S) \to \Spec(S_0)$ is a geometric quotient for periodic $S$. Thus the following notion of stability is well defined.} point of the action by $G$ associated to $d \in D_\BQ$, if $x \in X^\text{ss}(d)$, the orbit cone $\sigma(x)$ is maximal-dimensional and the orbit\footnote{For the definition of orbit, we refer to \cite{AT}, Section 3.1.4.} $Gx$ is closed in $X^\text{ss}(d)$. We will write $x \in X^\text{s}(d)$.

        \item An open subset $U \subseteq X^\text{gss} \subseteq \Spec(S)$ is said to be \emph{a stable subset of $X^\text{gss}$}, if 
        \begin{enumerate}[label=(\roman*)]
            \item  the action of $G$ on $U$ is proper,
            \item  the orbit of every $x \in U$ is maximal-dimensional (and closed), and
            \item $U$ is maximal with respect to inclusion among open subsets $V \subseteq X^\text{gss}$ satisfying (i) and (ii).
        \end{enumerate}
    \end{enumerate}
\end{definition}

\begin{example}\label{ex:double_origin}
  Let $S = \BC[x, y, z, w]$ and $D = \BZ^2$, where $\deg(x) = \deg(y) = (1, 0)$, $\deg(z) = (1, 1)$, $\deg(w) = (0, 1)$ and $\Gen^D(S) = \{xw, yw, zw, xz, yz\}$. We compute that $G= \BG_m^2$ acts on $\BA^4$ via
    \begin{align*}
        ((\lambda, \mu) (x, y, z, w)) \mapsto (\lambda x, \lambda y, \lambda \mu z, \mu w) ,
    \end{align*}
    giving the global coordinates 
    \begin{align*}
        \Proj^D(S) \equ \{(xw:yw:z) \mid \text{ not } x = y = 0, \text{ not } w=z=0\},
    \end{align*}
    i.e.\ we have coordinates of a $\BP^2$ (for $w = 1$) and for $z = 0$ we have extra coordinates $(x:y:0)$ of a $\BP^1$, so $\Proj^D(S) \subseteq \BP^2 \times \BP^1$.
Note that all points $(x, y, z, 0) \in D(xz) \cup D(yz)$ and all points $(0, 0, z, w) \in D(wz)$ are mapped to $(0:0:1)$. Thus we have to remove one of the sets in order to receive a subset of stable points (or a maximal separated subset).
In fact $\BP^2 = D_+(w) = \Spec(S_{(xw)}) \cup \Spec(S_{(yw)}) \cup \Spec(S_{(zw)})$ does not contain the points coming from $(x, y, z, 0)$ and $\Bl_p(\BP^2) = D_+(x, y) = \Spec(S_{(xw)}) \cup \Spec(S_{(yw)}) \cup \Spec(S_{(xz)}) \cup \Spec(S_{(yz)})$ does not contain the points $(0, 0, z, w)$. 
\end{example}

\begin{remark}\label{rem:max_sep_sub_stable}
\begin{enumerate}
    \item The definition of semistable points highly depends on $d$: As $1 \in X^\text{ss}(0)$, clearly all points of $\Spec(S)$ are semistable (see also \cite{CLS}, Example 14.1.7).

    \item The orbit cone $\sigma(x)$ is maximal-dimensional if and only if the stabilizer group scheme $\Stab_G(x)$ (which coincides with $G_x$) is finite (as scheme over $S_0$).

    \item We cannot define $X^s$ to be the union of all stable points $X^\text{s}(d)$ for $d \in D$, as in general we do not know if $G$ acts properly on the union. In fact, in Example~\ref{ex:double_origin} we already saw that it does not. In particular, we see that $D_+(x,y)$ and $D_+(w)$ each satisfy the condition of a stable subset. Hence, we have two distinct but maximal stable open subsets.

    \item Examining the affine projection morphism \begin{align*}
    \pi_+ \colon \Spec(S) \setminus V(S_+) \, \to\, \Proj^D(S) 
\end{align*}
    and Definition~\ref{def:semi_stable} (7) yields that each choice of a stable subset $U \subseteq X^\text{gss}$ is equivalent to the choice of a maximal separated subset of $\Proj^D(S)$:
    By (i), the topological space $U / G$ is Hausdorff. Then by (ii) (and \cite{paper1}, Proposition 2.15), $\pi_+(U)$ is a geometric quotient, identifying $U/G $ with $\pi_+(U)$. Consequently, $\pi_+(U)$ is a separated subset of $\Proj^D(S)$, as it is identified with a Hausdorff space. Finally, condition (iii) ensures that this separated subset is maximal with respect to inclusion.
\end{enumerate}
\end{remark}

The following Lemma is \cite{CRB}, Lemma 3.1.2.7. Since there is no proof in \cite{CRB}, we will provide one.

\begin{lemma}\label{lem:X^ss}
    Let $S$ be a multigraded ring.
    For each $d \in \sigma(S) \cap D_\BQ$, the set of semistable points is given by
    \begin{align*}
        X^{\text{ss}}(d) \deq \{x \in X \mid d \in \sigma(x)\}.
    \end{align*}
\end{lemma}

\begin{proof}
  Let $d \in \sigma(S) \cap D_\BQ$.
  For \LI, let $x \in X^{\mathrm{ss}}(d)$. Then there exist $n > 0$ and a homogeneous element $f \in S_{nd}$ such that $f(x) \neq 0$. Hence $nd \in \sigma(x)$ by the definition of the orbit cone $\sigma(x)$, and therefore also $d \in \sigma(x)$ since $\sigma(x)$ is a cone. 
  
  For \RI, let $x \in X$ such that $d \in \sigma(x)$. By definition of $\sigma(x)$, there exist homogeneous elements $f_1,\ldots,f_m \in S$ with $f_i(x)\neq 0$ and non-negative rational numbers $a_1,\ldots,a_m$ such that \[ d = \sum_{i=1}^m a_i \deg(f_i). \] Choose $N>0$ such that $Na_i \in \mathbb{Z}_{\geq 0}$ for every $i$. Then \[ Nd = \sum_{i=1}^m (Na_i)\deg(f_i) = \deg\left(\prod_{i=1}^m f_i^{Na_i}\right). \] Since each $f_i(x)$ is nonzero, we also have \[ \left(\prod_{i=1}^m f_i^{Na_i}\right)(x) \neq 0. \] Thus \[ \prod_{i=1}^m f_i^{Na_i} \in S_{Nd} \] is a homogeneous element not vanishing at $x$. Hence $x\in X^{\mathrm{ss}}(d)$.
\end{proof}

If we adopt this Lemma to our notion of geometric semistability, we can deduce the following:

\begin{corollary}\label{cor:gss_direct_lemma}
    Let $S$ be a $D$-graded ring, $x \in \Spec(S)$ and $d \in \sigma(S) \cap D_\BQ$. Then it holds that 
    \begin{align*}
        X^\text{gss}(d) \equ \{ x \in \Spec(S) \mid d \in \sigma(x) \text{ and } \dim(\sigma(x)) = \rank(D)\} \subseteq \Spec(S).
    \end{align*}    
\end{corollary}

Next, we show that the geometrically semistable locus is determined by the relevant elements.

\begin{corollary}\label{cor:X^ss_via_CRB}
    Let $S$ be a $D$-graded ring, $x \in \Spec(S)$ and $d \in \sigma(S) \cap D_\BQ$. 
\begin{enumerate}[label=(\roman*)]
    \item It holds
    \begin{align*}
        X^\text{gss}(d) \equ \{ x \in \Spec(S) \mid \exists n > 0, f \in \Rel^D(S) \cap S_{nd}: f(x) \neq 0\} \subseteq \Spec(S),
    \end{align*}    
    so that $x$ is geometrically semistable if and only if $x$ is contained in $D(f)$ for some relevant $f$.

    \item In particular, if $\sigma(x)$ is maximal-dimensional, then it can be written as the union of all relevant cones $\CC_D(f)$, where $x \in D(f)$, i.e.
    \begin{align*}
          \sigma(x) \equ \overline{\bigcup_{f \in \Rel^D(S), x \in D(f)} \CC_D(f)}.
    \end{align*} 
\end{enumerate}      
\end{corollary}

\begin{proof}
\begin{enumerate}[label=(\roman*)]
\item    \RI Let $x$ be contained in the right-hand side. Then there exists some relevant $f$ in some $S_{nd}$ such that $f(x) \neq 0$.  In particular, this holds for all divisors of $f$ and thus $\CC_D(f) \subseteq \sigma(x)$. By \cite{paper1}, Lemma 1.14 we know that $\CC_D(f)$ has maximal dimension, hence also $\sigma(x)$.

    \LI Conversely, let $x \in X^\text{gss}$. Then $d \in \sigma(x)$ and $\sigma(x)$ is maximal-dimensional. Thus, there are at least $r = \rank(D)$ many homogeneous elements $h_i$ satisfying $h_i(x) \neq 0$.
    As $\sigma(x)$ has full dimension, we can assume without loss of generality that the degrees $\deg(h_i)$ are pairwise linearly independent.
    Therefore $f := \prod_{i=1}^r h_i$ is relevant and does not vanish at $x$, hence $x \in X^\text{gss}(d)$.

\item Regarding the decomposition of $\sigma(x)$:

\LI Let $d \in \sigma(x) \cap D_\BQ$, where $\sigma(x)$ is of maximal-dimensional. Then there exists a relevant element $f := \prod_{i=1}^r h_i$ like above. 
As $d \in \sigma(x)$, one of the $\deg(h_j)$ (for some $j$) 
must lie in the same cone as $d$. In particular, we may arrange that $\deg(h_j)=d$. Now as $h_j \mid f$, it follows that
\begin{align*}
    d \in \CC_D(h_j) \subseteq \CC_D(f) \subseteq \overline{\bigcup_{f \in \Rel^D(S), f(x) \neq 0} \CC_D(f)}.
\end{align*}

\RI Let $d \in D$ be contained in the right-hand side, i.e.\ there exists some relevant $f \in S$ such that $x$ is contained in $D(f)$ and $d \in \CC_D(f)$. Hence, there exists some $N >0$ such that $Nd \in D^f$ and thus there exists some divisor $g$ of $f$ of degree $Nd$ such that $g(x) \neq 0$. We deduce that $x \in X^\text{ss}(d)$ and therefore $d \in \sigma(x)$.
\end{enumerate}
\end{proof}

\begin{remark}\label{rem:gss_and_gs}
    The proof of Corollary~\ref{cor:X^ss_via_CRB} has shown that every maximal-dimensional orbit cone $\sigma(x)$ contains the degree of a relevant element (and vice versa). Therefore, an element $x \in \Spec(S)$ is contained in a maximal-dimensional orbit cone $\sigma(x)$ if and only if there exists some relevant element $f$ satisfying $f(x) \neq 0$.
\end{remark}

If $\sigma(x)$ is polyhedral, we can simplify the previous corollary.

\begin{corollary}\label{cor:X^ss_via_CRB_finite}
Let $S$ be effectively $D$-graded, and $x \in \Spec(S)$. Assume that $\sigma(x)$ is maximal-dimensional and polyhedral.
Then there exist some relevant $f \in S$ such that $\sigma(x) = \CC_D(f)$.
\end{corollary}

\begin{proof}
    Let $\sigma(x)$ be polyhedral and maximal-dimensional, $r = \rank(D)$. Then there exist generators $d_1, \ldots, d_m$ of $\sigma(x)$, $m \ge r$, such that for each $d_i$ there exists a homogeneous element $h_i \in S$ of degree $d_i$ satisfying $h_i(x) \neq 0$. Then $f = \prod_{i=1}^m h_i$ is relevant. We claim that $\sigma(x) = \CC_D(f)$. By construction, we see that all $d_i$ are contained in $\CC_D(f)$, hence \LI holds. Conversely, let $d \in \CC_D(f) \cap D$. Then $d$ is the degree of some homogeneous element $g$ dividing some $f^k$ for $k \in \BN$. As $f^k(x) \neq 0$, also $g(x) \neq 0$ and thus $d \in \sigma(x)$.
\end{proof}

Altogether, we can state the following Theorem.

\begin{theorem}\label{thm:X^gss_as_relevant_cone}
    Let $D$ be a finitely generated abelian group, $S$ a $D$-graded ring and $G = \Spec(S_0[D])$. 
     The set of geometrically semistable points with respect to the action of the diagonalizable 
$S_0$-group scheme $G$ on $\Spec(S)$ is  given by
    \begin{align*}
        X^\text{gss} \equ \bigcup_{f \in \Rel^D(S)} D(f) \equ D(S_+).
    \end{align*}

    In particular, $\Proj^D(S)$ is characterized in terms of the geometric quotient $\Proj^D(S) := X^\text{gss} // G$, so it can be described without the notion of relevance.
\end{theorem}

\begin{proof}
    By Corollary~\ref{cor:X^ss_via_CRB} we have that $$\sigma(x) = \bigcup_{f \text{relevant}, f(x) \neq 0} \CC_D(f),$$ where we know that the relevant elements $f$ lie in some graded part $S_{nd}$ for $d \in \sigma(x)$ and $n > 0$. Now  as $X^\text{gss}$ is the union of all $X^\text{gss}(d)$ for $d \in D$, it holds
   \begin{align*}
       X^\text{gss} &\equ \bigcup_{d \in D} X^\text{gss}(d) \equ \{x \in \Spec(S) \mid \dim(\sigma(x)) = \rank(D)\} \\
       &\stackrel{\ref{cor:X^ss_via_CRB}}{\equ} \{x \in \Spec(S) \mid \exists f \in \Rel^D(S): f(x) \neq 0\}        \equ D(S_+).
   \end{align*} 
\end{proof}

We can also describe the GIT cones purely in terms of relevant elements (cf.\  \cite{CRB}, Lemma 3.1.2.10).

\begin{corollary}\label{cor:GIT_relevant_cones}
    The GIT cone for $d \in \lint(\sigma(S)) \cap D$ is given by
    \begin{align*}
        \lambda(d) \equ \bigcap_{d \in \lint(\sigma(x))} \sigma(x) \equ \bigcap_{\substack{f \in \Rel^D(S)  \\ d \in \lint(\CC_D(f))}} \CC_D(f).
    \end{align*}
    In particular, if $S$ is noetherian, integral and $S_+ = (f_1, \ldots, f_n)$ is radical, we can simplify
        \begin{align*}
        \lambda(d)   \equ \bigcap_{\substack{i \in \{1,\ldots, n\}  \\ d \in \lint(\CC_D(f_i))}} \CC_D(f_i).
    \end{align*}
 Likewise, the same holds true for the respective interiors.
\end{corollary}

\begin{proof}
The first equation holds by \cite{CRB}, Lemma 3.1.2.10. By Corollary~\ref{cor:X^ss_via_CRB}, we know that we can cover each maximal-dimensional orbit cone $\sigma(x)$ with the set of relevant elements not vanishing at $x$.

Take some relevant $f \in S$ with $d \in \lint(\CC_D(f))$, i.e.\ there exists a homogeneous divisor $g \in S_d$ of some $f^k$ with $g \in \sigma(x)$, as $f(x) \neq 0$. Hence $d \in \CC_D(g) \subseteq \CC_D(f)$.
As $f$ was arbitrary, we deduce $d \in \CC_D(f)$ for all $f \in S$ such that $d \in \lint(\sigma(x))$.

Conversely, let $e \in D \cap \lint(\CC_D(f))$ satisfy $e \in \lint(\CC_D(f))$ for all relevant $f$ such that $d \in \lint(\CC_D(f))$. For any $x$ such that $d \in \lint(\sigma(x))$, there exists some relevant $h \in S$ with $h(x) \neq 0$ and $d \in \lint(\CC_D(h))$. Thus $e \in \CC_D(h) \subseteq \sigma(x)$. But $x$ was arbitrary, so $e$ has to be contained in the intersection of all those $\sigma(x)$, that is, $e \in  \bigcap_{d \in \lint(\sigma(x))} \sigma(x)$.

Now let $S$ be noetherian and integral such that $S_+ = (f_1, \ldots, f_n)$ is radical. As $d \in \lint(\sigma(S)) \cap D$, every orbit cone living in the intersection is maximal-dimensional and polyhedral (by assumption). It's also clear that each $f_i$ defines a maximal-dimensional orbit cone $\sigma_i(x)$. Now let $\sigma(x)$ be a maximal-dimensional orbit cone, given by $\CC_D(f)$ for some relevant $f$. Since $f \in S_+ = (f_1,\ldots, f_n)$ and $S_+$ is radical, it holds $D(f) = \bigcup_{i=1}^n D(f_i)$. As we assumed $S$ to be an integral domain, $D(f)$ is irreducible, and hence $D(f) = D(f_j)$ for some $j \in \{1,\ldots, n\}$. But this implies that $\CC_D(f) = \sigma_j(x)$. In particular,
\begin{align*}
    \lambda(d)   \equ \bigcap_{i:d \in \lint(\sigma(x))} \sigma(x) \equ \bigcap_{d \in \lint(\sigma_i(x))} \sigma_i(x) \equ  \bigcap_{\substack{i :d \in \lint(\CC_D(f_i))}} \CC_D(f_i).
\end{align*}
\end{proof}

\begin{remark}\label{rem:arbitrary_GIT_cone}
    In the sense of Definition~\ref{def:semi_stable} (2), we might call $\lambda(d)$ a \emph{geometric GIT cone} if 
        \begin{align*}
            \lambda(d) \equ \bigcap_{\substack{x\in\Spec(S), d \in \sigma(x) \\ \sigma(x) \text{ max. dim.}}} \sigma(x).
        \end{align*}
We can also describe arbitrary GIT cones: For $d \in \sigma(S)$ we define
\begin{align*}
    \lambda(d) \equ \bigcap_{\substack{g \in S \text{ hom.} \\ d \in \CC_D(g)}} \CC_D(g).
\end{align*}
\end{remark}

\begin{remark}
 Note that by \cite{CRB}, Theorem 3.1.2.8, the GIT cones $\lambda(d)$ form a quasifan (in $D_\BQ$) if $S$ is a finitely generated $S_0$-algebra for $S_0$ an algebraically closed field of characteristic zero.
 Therefore, if $X_\Sigma$ is a complex quasiprojective toric variety such that $\Sigma$ has full-dimensional convex support, the above-mentioned corresponding quasifan coincides with the secondary fan (after scalar extension).
    In particular, the chambers of the secondary fan are precisely the GIT cones and hence are fully determined by relevant elements.
    We will discuss this in more detail in Section~\ref{sec:VGIT}.
Note that the same holds true for $X$ being a Mori dream space. In this case, the quasifan is precisely the Mori chamber decomposition of $X$ (see \ref{sec:MDS} for more details).
\end{remark}

\begin{example}\label{ex:standard_ex_stable_points}
Let $S = \BC[x, y, z]$ be graded by $\BZ\times \BZ/2\BZ$, such that $\deg(x) = (1, 0)$, $\deg(y) = (0, 1)$,  $\deg(z) = (1, 1)$, hence $\Gen^D(S) = \{x, z\}$.
    We have an action $\BG_m \times \mu_2 \looparrowright \BA^3$, given by
    \begin{align*}
        (\lambda, \psi) (x, y, z) \mapsto (\lambda x, \psi y, \lambda \psi z) .
    \end{align*}
    Points of type $(0, y^2, 0)$ are $G$-invariant as well as points $(x^2, yz, z^2)$ and points $(x, yz, 0)$ as well as $(x, 0, yz)$ and $(xy, yx, z)$. 
    Following \cite{paper1}, Theorem A, we compute $S^{\mu_2} = \BC[x, y^2, z^2, yz]$. Then the induced $\BG_m$ action on $\Spec(S^{\mu_2})$ is given by $t \cdot x = tx$, $t \cdot y^2 = t^2y^2$, $t \cdot z^2 = t^2z^2$ and $t \cdot yz = t yz$. Hence, we can identify
    \begin{align*}
        \Proj^D(S) \equ \Proj^\BN(\BC[a, b, c, d]/(d^2-bc)),
    \end{align*}
    via $a \mapsto x$, $b\mapsto y^2$, $c \mapsto z^2$ and $d \mapsto yz$, where $\deg(a) = 1, \deg(b) = 0, \deg(c) = 2$ and $\deg(d) = 1$.
    As $\Proj^D(S)$ is separated, it has a unique open subset of stable points, which is $\Proj^D(S)$ itself.
    But note that, even though $\Proj^D(S)$ is separated, $\pi_+$ is not a Zariski locally trivial pseudo $G$-torsor, as it is for Example~\ref{ex:double_origin} (cf.\ \cite{paper1}, Theorem B).
\end{example}




\section{Toric VGIT via $D$-graded Proj}\label{sec:VGIT}
Given a quasiprojective normal toric variety $X_\Sigma$ such that $\Sigma$ has full-dimensional convex support (over an algebraically closed field $k$ of characteristic zero), \cite{GKZ} studied actions of closed subgroups $H \le \Spec(\Cox(X)_0[\Cl(X)]) \cong \BG_{m, \Cox(X)_0}^n$ (for $n = \rank(\Cl(X)) + \dim(X)$). Using a character $\chi$ of $G$, one can lift the action of $H$ to a trivial line bundle over $\Spec(\Cox(X_\Sigma))$.
Full dimensional convex support means that the support of $\Sigma$ is a convex subset of $N_\BR$ satisfying $\dim(|\Sigma|) = \dim(N_\BR)$ (see \cite{CLS}, page 265).

The information on the different GIT quotients arising from different characters $\chi$ is stored in the so-called \emph{secondary fan} of $X_\Sigma$. There are many different ways to think about the secondary fan, as it is shown in \cite{CLS}, Chapters 14 and 15. We will give a quick introduction to the construction of the secondary fan, following \cite{CLS} Chapter 14. Then we will show how to construct the secondary fan and its chambers purely in terms of multigraded objects, if $\Sigma$ is simplicial.

Concretely, we can associate with each chamber in the secondary fan an irrelevant ideal, which is given in terms of the relevant cones $\CC_D(f)$ containing any element of the interior of the chamber.
Using this fact, we show that $\Proj^{\Cl(X_\Sigma)}(\Cox(X_\Sigma))$ is given by the direct limit of the different GIT models being classified by the secondary fan.
In particular, this result gives an answer to the question of whether the Proj construction is an inverse construction to the Cox construction.

Finally, we generalize the construction of the secondary fan and the movable cone to simplicial toric prevarieties.
We follow \cite{CLS} closely.
For this chapter, $S$ will always be a polynomial ring with finitely many variables over $k$. Additionally, toric varieties are assumed to be normal.

\subsection{Toric GIT}\label{sec:toric_GIT}

Let $X = X_\Sigma$ be a quasiprojective toric variety such that $\Sigma$ has full-dimensional convex support. Then
\begin{align*}
    \widehat{G} \deq \{\chi\colon G\to \BG_m \mid \chi \text{ is a homomorphism of algebraic groups }\}
\end{align*}
denotes the character group of $G := \Spec(\Cox(X)_0[\Cl(X)])$. There is a natural map $\gamma\colon \BZ^n \to \widehat{G}, a \mapsto \chi^a$. Hence for $M = \ker(\gamma)$, we get a short exact sequence
\begin{align}\label{eq:exact_sequence}
    0 \ \to \ M \ \stackrel{\delta}{\to} \ \BZ^n \ \stackrel{\gamma}{\to} \ \widehat{G} \ \to 0
\end{align}
(cf.\  \cite{CLS}, Lemma 14.2.1 (a)). The images of the standard basis elements $e_i \in \BZ^r$ will be denoted by $\beta_i$, i.e.\
\begin{align*}
    \beta_i \deq \gamma(e_i), \ \ i = 1,\ldots, n.
\end{align*}
The dual of $\delta$ is the map $\delta^\ast\colon\BZ^n \to N$, where $N = M^\ast$ is the dual of $M$. The images of the standard basis elements $e_i$ will be denoted by $\nu_i$, i.e.\
\begin{align*}
    \nu_i \deq \delta^\ast(e_i), \ \ i = 1,\ldots, n ,
\end{align*}
and hence
\begin{align*}
    \delta_\BR(m) \equ (\langle m, \nu_1 \rangle,  \ldots \langle m, \nu_n \rangle ), \ m \in M_\BR.
\end{align*}
We denote the dual short exact sequence by
\begin{align*}
    0 \ \to \ \widehat{G}^\ast \ \stackrel{\gamma^\ast}{\to} \ \BZ^n \ \stackrel{\delta^\ast}{\to} \ N \ \to 0,
\end{align*}
where 
\begin{align*}
    (\gamma^\ast)_\BR(v) = (\langle \beta_1, v \rangle, \ldots, \beta_n, v\rangle), v \in \widehat{G}^\ast
\end{align*}
by duality.
The two vector configurations $\nu_1, \ldots, \nu_n$ and $\beta_1, \ldots, \beta_n$ are said to be \emph{Gale dual} to each other. For simplicity, we will denote $C_\beta := \Cone(\beta_1,\ldots, \beta_n) \subseteq \widehat{G}_\BR$ and $C_\nu := \Cone(\nu_1,\ldots, \nu_n) \subseteq N_\BR$ as in \cite{CLS}.

The core statement from toric GIT is \cite{Cox}, Theorem 2.1 (see Theorem~\ref{thm:cox_quot} below). Let 
\begin{align*}
    \Cox(X)_+^C \deq \left( \prod_{\rho \not\in\sigma(1)} x_\rho \mid \sigma \in \Sigma \right)
\end{align*}
 be the irrelevant ideal as defined by \cite{Cox}, where $\sigma(1)$ denotes the set of one-dimensional faces of $\sigma \in\Sigma$. 

\begin{theorem}\label{thm:cox_quot}
    Let $X = X_\Sigma$ be a toric variety such that $\Sigma$ has full-dimensional convex support, and $S = \Cox(X)$. Then it holds:
    \begin{enumerate}[label=(\roman*)]
        \item $\Spec(S) \setminus V(S_+^C)$ is invariant under the action of $G$.
        \item $X$ is naturally isomorphic to the categorical quotient of $\Spec(S) \setminus V(S_+^C)$ by $G$.
        \item The quotient from (ii) is geometric if and only if $\Sigma$ is simplicial.
    \end{enumerate}
\end{theorem}

So the idea is to take $\Spec(\Cox(X))$, remove some points, and then mod out by $G$. However, instead of looking at the action of all of $G$ at once, we can take a random `direction' $d \in \Cl(X)$ and then compute the points that have to be removed along this direction.
In practice one lifts the $G$-action to the rank $1$ trivial vector bundle $\BA^n \times \BA \to \BA^n$ associated to $\chi$, where $g \in G$ acts on $(p, t) \in \BA^n \times \BA$ via
\begin{align}\label{eq:linearized_character}
    g \cdot (p, t) \ \equ \ (g \cdot p, \chi(g)t).
\end{align}
The corresponding sheaf of sections is denoted by $\CL_\chi$ and is called the \emph{linearized line bundle} with character $\chi$.
Note that $f \in \Gamma(\BA^n, \CL_\chi)^G$ if and only if $f$ satisfies Equation~\ref{eq:1} (for more details, we refer to \cite{CLS}, Lemma 14.1.1 and \cite{paper1}, Appendix A).

For simplicity, we will use the notation from \cite{CLS}.

\begin{definition}\label{def:(semi)stable_points_toric}
Let $G \subseteq \BG_m^n$ and $\chi \in \widehat{G}$.
The set of semistable (resp. stable) points is denoted by $(\BA^n)_\chi^{\text{ss}}$ (resp $(\BA^n)_\chi^{\text{s}}$).
\end{definition}

Now \cite{CLS} use the $\BN$-graded ring 
\begin{align*}
    R_\chi \deq \bigoplus_{d= 0}^\infty \Gamma(\BA^n, \CL_{\chi^d})^G,
\end{align*}
that is a finitely generated $k$-algebra by \cite{CLS}, Lemma 14.1.10, to define the GIT quotient associated to a character $\chi$:

\begin{definition}\label{def:GIT_quot_via_N_Proj}
    The GIT quotient of $\chi \in \widehat{G}$ is defined by
    \begin{align*}
        \BA^n \sslash_\chi G \equ \Proj^\BN(R_\chi).
    \end{align*}
\end{definition}

Hence, VGIT means that we study the GIT quotients arising from various characters (see \cite{T} and \cite{DH} for more details on this subject).
In the following sections, \cite{CLS} then give several equivalent ways to define the above GIT quotient. We are mostly interested in the approach via irrelevant ideals, hence we will only give a quick overview of the theory.

\subsection{Secondary Fan}\label{sec:sec_fan}
We assume that $X_\Sigma$ is a quasiprojective simplicial normal toric variety (over $k$) such that $\Sigma$ has full-dimensional convex support.
We try to outline the objects that play a key role in the secondary fan framework. 

Note that we restrict to the simplicial case due to Theorem~\ref{thm:cox_quot} (iii) and \cite{paper1}, Theorem A (here, also \cite{BS}, Lemma 2.1 applies). Concretely, $\Proj^D(S)$ is a geometric quotient for any multigraded ring $S$.

\begin{definition}
For each element $(a_1, \ldots, a_n) = a \in \BZ^n$ we define:
\begin{enumerate}[label=(\arabic*)]
    \item a character $\chi = \chi^a$, where we denote by $\chi^a$ the character defined by $\chi^{\gamma(a)}$.
    \item a polyhedron $P_a := \{m \in M_\BR \mid \langle m, \nu_i \rangle \ge -a_i,\ 1 \le i \le n\} \subseteq M_\BR$, where for fixed $\chi$, the various choices of $a$ differ by elements of $\delta(M)$ and the toric variety associated to $P_a$ is exactly $\Proj(R_{\chi^a})$.
    \item a graded ring $k[C(P_a) \cap (M \times \BZ)] \cong R_{\chi^a}$, where $C(P_a)$ is the cone over $P_a$ (cf.\  \cite{CLS}, Theorem 14.2.13).
    \item a polyhedron $P_\chi := \{b \in \BR^n_{\ge 0} \mid \gamma_\BR(b) = \chi \otimes 1\} = \gamma_\BR^{-1}(\chi \otimes 1) \cap \BR^n_{\ge 0}$, where $P_\chi = \delta_\BR(P_a) + a$ depends only on the class of $a$ as element in $\Cl(X)_\BR$.
\end{enumerate}
\end{definition}

While $P_a$ is not a lattice polyhedron in general and hence cannot give rise to the GIT quotient $\BA^n \sslash_{\chi^a} G$ (cf.\  \cite{CLS}, Example 14.2.3), it holds that $\dim(\BA^n \sslash_\chi G) = \dim(P_\chi)$ (cf.\  \cite{CLS}, Corollary 14.2.16).

\begin{definition}\label{def:pair_GKZ}
    For $a \in \BZ^n$ and $1 \le i \le n$, we call the set
    \begin{align*}
        F_{i,a} \deq \{m \in M_\BR \mid \langle m, \nu_i \rangle = -a_i\}
    \end{align*}
    a \emph{virtual facet} of $P_a$. Consequently, the virtual facets of $P_\chi$ are defined by
    \begin{align*}
        F_{i, \chi} \deq P_\chi \cap V(x_i),
    \end{align*}
    where we use that the Cox ring of a toric variety is a polynomial ring with $n = \rank(\Cl(X))+ \dim(X)$ many variables (cf.\  \cite{HK}, Corollary 2.10). The empty virtual facets of $P_a$ are determined by the set
    \begin{align*}
        I_\emptyset \deq \{i \in \{1, \ldots, n\} \mid F_{i,a} \equ \emptyset \}.
    \end{align*}
    Thus, $a$ gives rise to a pair $(\Sigma, I_\emptyset)$, where $\Sigma$ is the normal fan of $P_a$ in $N_\BR$, such that
    \begin{enumerate}[label=(\alph*)]
        \item $|\Sigma| = C_\nu$.
        \item $X_\Sigma$ is a quasiprojective toric variety such that $\Sigma$ has full-dimensional convex support.
        \item $\sigma = \Cone(\nu_i \mid \nu_i \in \sigma, i \not\in I_\emptyset)$ for every $\sigma \in \Sigma$.
    \end{enumerate}
\end{definition}

The virtual facets determine the semistable points and the irrelevant ideal of $\chi^a$:

\begin{proposition}
    Let $p = (p_1,\ldots, p_n) \in k^n$ and set $I(p):=\{i \mid 1 \le i \le n, p_i =0\}$. Then $p \in (\BA^n)^\text{ss}_\chi$ if and only if $\cap_{i \in I(p)} F_{i,\chi} \neq \oldemptyset$.
\end{proposition}

\begin{proof}
    By \cite{CLS}, Proposition 14.2.21.
\end{proof}
An algebraic interpretation of the previous result gives the irrelevant ideal associated to $\chi = \chi^a$.

\begin{definition}\label{def:toric_irrel_ideal}
    Let $\chi \in \widehat{G}$. We define the \emph{irrelevant ideal} of $\chi$ to be 
    \begin{align*}
        B(\chi) \deq \left(\prod_{i \not\in I} x_i \mid I \subseteq \{1, \ldots, n\},\ \bigcap_{i\in I} F_{i, \chi} \neq \oldemptyset\right) \ \subseteq \ S = k[x_1, \ldots, x_n] .
    \end{align*}
    The vanishing locus of $B(\chi)$ is the \emph{exceptional set} $Z(\chi) = V(B(\chi))$.
\end{definition}

\begin{corollary}\label{cor:N_Proj_via_ss}
    For $\chi \in \widehat{G}$ it holds $(\BA^n)_\chi^\text{ss} = \BA^n\setminus Z(\chi)$, so that
    \begin{align*}
        \BA^n \sslash_\chi G \stackrel{\sim}{\equ} (\BA^n \setminus Z(\chi)) / G.
    \end{align*}
\end{corollary}

\begin{proof}
    By \cite{CLS}, Exercise 14.2.9.
\end{proof}

Viewing $\chi$ as element of $\widehat{G}_\BR$ via $\chi \otimes 1$ give rise to the following classifications:

\begin{proposition}
    Let $\chi \in \widehat{G}$. Then $(\BA^n)_\chi^\text{ss} \neq \oldemptyset$ if and only if $\chi \otimes 1 \in C_\beta$. 
    If $\chi \otimes 1 \in C_\beta$, then $(\BA^n)_\chi^\text{s} \neq \oldemptyset$ if and only if $\chi \otimes 1 \in \lint(C_\beta)$.
\end{proposition}

\begin{proof}
    By \cite{CLS}, Proposition 14.3.5 and 14.3.6.
\end{proof}

\begin{remark}
    If $\chi \otimes 1 \in C_\beta$, then $\chi \otimes 1 \in \lint(C_\beta)$ is also equivalent to $F_{i, \chi}$ being a proper face of $P_\chi$ for all $i$. In particular, $(\BA^n)_\chi^\text{s} \neq \oldemptyset$ if and only if $\BA^n \sslash_\chi G$ has the correct dimension, with is $n - \dim(G)$ (cf.\ \cite{paper1}, Lemma 2.20).
\end{remark}

The central notion in the context of the secondary fan of a toric variety is that of a \emph{generic} character (cf.\  \cite{CLS}, Definition 14.3.14).

\begin{definition}
    A character $\chi \in \widehat{G}$ is called \emph{generic}, if $\chi \otimes 1 \in C_\beta$ and for every subset $\beta'\subseteq \{\beta_1, \ldots, \beta_n\}$ with $\dim(\Cone(\beta')) < \dim(C_\beta) = \dim(G)$, we have $\chi \otimes 1 \not\in\Cone(\beta')$.
\end{definition}

Generic characters are characterized as follows.

\begin{theorem}\label{thm:generic_char_equiv}
    For $\chi \otimes 1 \in C_\beta$, the following are equivalent:
\begin{enumerate}[label=(\alph*)]
    \item $\chi$ is generic.
    \item Every vertex of $P_\chi$ has precisely $\dim(G)$ nonzero coordinates.
    \item $P_\chi$ is simple of dimension $n-\dim(G)$, every virtual facet $F_{i,\chi} \subseteq P_\chi$ is either empty or a genuine facet, and $F_{i,\chi} \neq F_{j,\chi}$ if $i \neq j$ and $F_{i, \chi}, F_{j, \chi}$ are nonempty.
    \item $(\BA^n)_\chi^\text{s} \equ (\BA^n)_\chi^\text{ss}$.
\end{enumerate}
\end{theorem}

\begin{proof}
    By \cite{CLS}, Theorem 14.3.14.
\end{proof}

\begin{remark}\label{rem:X^gss_already_proof}
    Note that condition (d) in the previous theorem already implies, that generic characters $\chi^a$ define a subset $(\BA^n)_\chi^\text{ss} \subseteq X^\text{gss}$, as the equation $(\BA^n)_\chi^\text{s} \equ (\BA^n)_\chi^\text{ss}$ forces the orbit cones to be maximal dimensional, so to be given in terms of relevant elements by Corollary~\ref{cor:X^ss_via_CRB_finite}.
\end{remark}

The bottleneck here is to show that $C_\beta$ has the structure of a fan. The crucial part is to construct the cones of $C_\beta$. This is done in \cite{CLS}, Chapter 14.4, via the so-called \emph{GKZ cones} (Gelfand--Kapranov--Zelevinsky).
Such a cone is defined for a pair $(\Sigma, I_\emptyset)$ from Definition~\ref{def:pair_GKZ} in terms of convex support functions on $\Sigma$, denoted $\text{CSF}(\Sigma)$:
\begin{align*}
\Tilde{\Gamma}_{\Sigma,I_\emptyset}
:=\left\{(a_1,\ldots,a_n)\in\mathbb{R}^n \;\middle|\;
\begin{aligned}
  &\exists\varphi\in\mathrm{CSF}(\Sigma):\ \varphi(\nu_i)=-a_i \quad (i\notin I_\emptyset),\\
  &\text{and }\ \varphi(\nu_i)\ge -a_i \quad (i\in I_\emptyset)
\end{aligned}
\right\}.
\end{align*}
Note that since we assume $X$ to have a fan with full-dimensional convex support, we know that all convex support functions on $X$ that are associated to a Cartier divisor are in fact strictly convex (cf.\ \cite{CLS}, Definition 6.1.12).
Now the actual \emph{GKZ cones} of the secondary fan are given by the cones 
\begin{align*}
    \Gamma_{\Sigma, I_\emptyset} \equ \Tilde{\Gamma}_{\Sigma, I_\emptyset} / \ker(\gamma_\BR) \subseteq \widehat{G}_\BR,
\end{align*}
when $\Sigma$ is a generalized fan. We refer to \cite{GKZ}, \cite{OP}, and \cite{HKP} for more details on this subject.
However, the chambers of the secondary fan can also be described using the irrelevant ideals only (cf.\  \cite{CLS}, Corollary 14.4.15). As this approach does not require the existence of support functions, we will only work with GKZ cones if necessary.

The most important fact to know is that $\chi \otimes 1$ lies in a maximal cone of $C_\beta$, which we call a \emph{chamber} of the secondary fan, if and only if $\chi$ is generic (cf.\  \cite{CLS}, Proposition 14.4.9). Thus, generic characters contain all the necessary information to describe the VGIT of $X$. In particular, there can only be finitely many distinct generic characters (modulo isomorphisms), as all $B(\chi)$ are combinatorial objects and there are only finitely many admissible choices. In fact, we will see that the Brenner--Schröer irrelevant ideal $S_+$ contains all irrelevant ideals $B(\chi)$ for generic $\chi$. In particular, the following characterization for generic immediately shows the presence of relevant elements.

In our setting, for a polynomial ring with finitely many variables $S = k[x_1,\ldots, x_n]$, we call a relevant element $f \in S$ \emph{monomic relevant} if there exists a subset $I \subseteq \{1, \ldots, n\}$ of size $|I|= \rank(D)$ such that 
\begin{align*}
    f \equ \prod_{i\in I} x_i.
\end{align*}
In \cite{paper1}, Lemma 1.19, we showed that $S_+$ is generated by monomic relevant elements in the setting of this chapter.

\begin{proposition}\label{prop:generic_iff_gen_dim(G)}
    Let $\chi \otimes 1 \in C_\beta$. Then $\chi$ is generic of and only if every minimal generator of $B(\chi)$ is monomic relevant.
\end{proposition}

\begin{proof}
    By \cite{CLS}, Proposition 14.4.14 (c), the generators of $B(\chi)$ for generic $\chi$ are given by elements of \emph{degree} $\dim(G)$, which is the same as saying that generators have a length $\dim(G)$ factorization (cf.\ \cite{paper1}, Definition 1.17), as the Cox ring of a quasiprojective toric variety $X_\Sigma$ such that $\Sigma$ has full-dimensional convex support is a polynomial ring. In particular, a minimal generator of $B(\chi)$ is monomic relevant by the definition of $B(\chi)$ (cf.\ Definition~\ref{def:toric_irrel_ideal}).
\end{proof}


\begin{theorem}\label{thm:secondary_relevant}
Let $X_\Sigma$ be a quasiprojective simplicial normal toric variety such that $\Sigma$ has full-dimensional convex support, $S = \Cox(X)$, and let $\chi = \chi^a$ be a character. Then it holds:
\begin{enumerate}[label=(\alph*)]
    \item $\chi^a$ is generic if and only if
    \begin{align*}
        B(\chi^a) \equ ( f\in \Gen^D(S) \mid a \in \CC_D(f) ) \unlhd S.
    \end{align*}
    \item The chamber $\sigma_a$ in the secondary fan containing $a$ is given by
    \begin{align*}
           \sigma_a \equ \bigcap_{\substack{f \in \Gen^D(S) \\ a \in \CC_D(f)}} \CC_D(f).
    \end{align*}
\end{enumerate}    
\end{theorem}

\begin{proof}
\begin{enumerate}[label=(\alph*)]
\item \Ra Let $\chi = \chi^a$ be a generic character, where $\chi \otimes 1 \in \lint(C_\beta)$. Then $B(\chi^a)$ is generated by monomic relevant elements by Proposition~\ref{prop:generic_iff_gen_dim(G)}.  Let $f \in B(\chi)$ be a minimal generator, so that there exists a factorization $f =\prod_{i \not\in I} x_i$ for some subset $I \subseteq \{1, \ldots, n\}$ such that $\bigcap_{i\in I} F_{i, \chi} \neq \oldemptyset$ by Definition~\ref{def:toric_irrel_ideal}, where $|I| = r := \rank(\Cl(X))$. As $\chi$ is generic, we may assume without loss of generality that $a_i = 0$ for $i \not\in I$ by the proof of \cite{CLS}, Theorem 14.3.14. Hence,
\begin{align*}
    \chi^a \otimes 1 \equ \sum_{i \in I} a_i \beta_i \equ \prod_{i\in I} x_i^{a_i}.
\end{align*}
Now we can see that $a$ lies in the interior of the cone generated by the rays corresponding to $x_i$ for $i \in I$. But this means $a \in \CC_D(f)$ by definition.

\La Let the equation hold. Now, as $B(\chi^a)$ is generated by relevant elements and $S$ is a polynomial ring, we know that $B(\chi^a)$ as a subset of $S_+$ is generated by monomic relevant elements by \cite{paper1}, Lemma 1.19, and thus we can apply Proposition~\ref{prop:generic_iff_gen_dim(G)}.

\item By construction, the chambers in the secondary fan corresponding to $\chi = \chi^a$ coincide with the GIT cones $\lambda(a)$ from Definition~\ref{def:semi_stable} (2), as the first is just a special case of the latter. In particular, the claim follows by Corollary~\ref{cor:GIT_relevant_cones}.
\end{enumerate}
\end{proof}



\begin{example}\label{ex:standard_ex_chambers}
    Consider the toric variety $X = \Bl_p(\BP^2)$. One computes that
    $\Cox(X) = S = \BC[x, y, z, w]$ and $D = \BZ^2$, where $\deg(x) = \deg(y) = (1, 0)$, $\deg(z) = (1, 1)$ and $\deg(w) = (0, 1)$ (cf.\ \cite{CLS} Example 14.2.17 for $r = -1$). Then $\Gen^D(S) = \{xw, yw, zw, xz, yz\}$. As $D = \BZ^2$, computing the chambers is rather easy. We can immediately see that we have three boundary rays and only one ray in the interior of $C_\beta$, which is coming from $z$. So we have the two chambers $C_1 = \overline{\Cone}(e_1, e_1+e_2)$, which is the intersection of $\CC_D(f)$ for $f = xw, yw, xz, yz$, and $C_2 = \overline{\Cone}(e_2, e_1+e_2)$, obtained by intersecting $\CC_D(f)$ for $f = xw, yw, zw$. Hence, the chambers correspond to the stable subsets computed in Example~\ref{ex:double_origin}. In particular, we see that $\Proj^{\Cl(X)}(\Cox(X))$ is $X$ glued with $\BP^2$.
\end{example}

\begin{theorem}\label{thm:list_equiv_chamber_var}
    Let $X_\Sigma$ be a quasiprojective simplicial normal toric variety such that $\Sigma$ has full-dimensional convex support, $S = \Cox(X)$, and let $\chi = \chi^a$ be a generic character. Then the following objects coincide:
    \begin{enumerate}[label=(\roman*)]
        \item $X_\Sigma$, where $\Sigma$ is the normal fan of $P_a$.
        \item $\Proj^\BN(R_{\chi^a})$.
        \item $\Proj^{\Cl(X)}(\sigma_a \cap \Cl(X)) = \Proj^{\Cl(X)}(\bigoplus_{d \in \sigma_a \cap \Cl(X)} H^0(X, d))$.
        \item $\Proj^{\Cl(X)}(k[h_1, \ldots, h_r])$, where $\deg(h_i)$ are the ray generators of $\sigma_a$.
        \item $\Proj^{D}_B(S)$, for $B = B(\chi^a)$ and $D = \Cl(X)$.
    \end{enumerate}
\end{theorem}

\begin{proof}
    The equivalence of (i) and (ii) holds by \cite{CLS}, Theorem 14.2.13. Let $\chi = \chi^a$ be a generic character. Recall that $\Proj^\BN(R_{\chi^a}) = (\BA^n \setminus Z(\chi^a)) / G$ by Corollary~\ref{cor:N_Proj_via_ss}. Then $B(\chi^a) = (f\in \Gen^D(S) \mid a \in \CC_D(f))$ by Theorem~\ref{thm:secondary_relevant}. In particular, $\Proj^{\Cl(X)}(\sigma_a \cap \Cl(X))$ has to be given by $(\Spec(S) \setminus V(M_a)) / G$, where $M_a$ is the set of relevant elements corresponding to $B(\chi^a)$. Now by construction, $Z(\chi^a) = V(M_a)$ and the equivalence of (ii) and (iii) holds true.
The equivalence of (iii) and (iv) is immediate, as the relevant locus of the ring in (iv) is exactly the set of generators of $B(\chi^a)$ in (iii).
Since the irrelevant ideal of the multigraded ring in (iv) is exactly $B(\chi^a)$, the affine charts in (iv) and (v) coincide with the same gluing data. Therefore, the schemes are isomorphic. The claim follows.
\end{proof}

\subsection{Nef and Movable Cones}
Let $X_\Sigma$ be a quasiprojective simplicial normal toric variety such that $\Sigma$ has full-dimensional convex support.
We can describe the \emph{movable cone} in terms of relevant elements of certain length. The \emph{nef cone} however, determines the chamber in the secondary fan of $X_\Sigma$, which $D$-graded Proj recovers $X_\Sigma$.
The following definitions are from \cite{CLS}, §6.

\begin{definition}\label{def:N^1N_1}
    For a normal variety $X$, define
    \begin{align*}
        N^1(X) \equ (\text{CDiv}(X) / \equiv) \otimes_\BZ \BR \text{  and  } N_1(X) \equ (Z_1(X)/\equiv )\otimes_\BZ \BR,
    \end{align*}
    where $Z_1(X)$ is the free abelian group generated by irreducible complete curves $C \subseteq X$. Elements of $Z_1(X)$ are called \emph{proper $1$-cycles}. 
\end{definition}

There is a nondegenerate bilinear pairing $N^1(X) \times N_1(X) \to \BR$. In fact, both $N^1(X)$ and $N_1(X)$ are finite-dimensional vector spaces over $\BR$ that are naturally their respective dual vector spaces.
Moreover, a Cartier divisor $D$ on $X$ is called \emph{nef} if $D\cdot C \ge 0$ for every irreducible curve $C \subseteq X$, where $D \cdot C$ denotes the nondegenerate bilinear pairing.

Now for toric varieties $X_\Sigma$ such that $\Sigma$ has full-dimensional convex support, the nef cone is easy to understand, as nef and base point free are equivalent in this setting (cf.\ \cite{CLS}, Theorem 6.3.12).
For general background on this subject, we refer to \cites{PAG, KM}.

\begin{definition}\label{def:nef_mori_cone}
    Let $X$ be a normal variety.
    \begin{enumerate}[label=(\alph*)]
        \item $\Nef(X)$ is the cone in $N^1(X)$ generated by classes of nef Cartier divisors. We call $\Nef(X)$ the \emph{nef cone}.
        \item $\NE(X)$ is the cone in $N_1(X)$ generated  by classes of irreducible complete curves.
        \item The \emph{Mori cone} $\overline{\NE}(X)$ is defined to be the closure of $\NE(X)$ in $N_1(X)$.
    \end{enumerate}
\end{definition}

The nef and Mori cones are closed convex cones and dual to each other, for example by \cite{CLS}, Lemma 6.3.19. In particular, $\NE(X)$ has full dimension in $N_1(X)$ and $\Nef(X)$ is strongly convex in $N^1(X)$. If $X$ has a secondary fan, $\Nef(X)$ is a rational polyhedral cone in $\Pic(X)_\BR$ by \cite{CLS}, Theorem 6.3.20.
If $X$ is also projective, the nef cone fully determines ampleness, i.e.\ a Cartier divisor on $X$ is ample if and only if its class in $\Pic(X)$ lies in the interior of $\Nef(X)$ by \cite{CLS}, Theorem 6.3.22.
Note that by \cite{CLS}, Proposition 15.1.3 it holds
\begin{align*}
    \Gamma_{\Sigma,I_\emptyset} \stackrel{\sim}{\equ} \Nef(X_\Sigma) \times \BR_{\ge 0}^{I_\emptyset}.
\end{align*}
So if we start with a (complex) toric variety $X = X_\Sigma$, the unique chamber giving $X$ is given by $\Gamma_{\Sigma, \emptyset}$, if $(\Sigma, \emptyset)$ is the GKZ-pair associated to $\Sigma$ (see Definition~\ref{def:pair_GKZ}). But when does this happen? We say that $\nu = (\nu_1, \ldots, \nu_n)$ is \emph{geometric}, if all $\nu_i$ are nonzero and generate distinct rays in $\Cl(X)_\BR$.

\begin{proposition}\label{prop:geometric}
    The secondary fan has a chamber with $I_\emptyset = \emptyset$ if and only if $\nu$ is geometric.
\end{proposition}

\begin{proof}
    By \cite{CLS}, Proposition 15.1.6.
\end{proof}

If we take the cone of all such pairs, we get the so-called \emph{moving cone} of the secondary fan. 

\begin{definition}\label{def:moving_GKZ_cone}
    Let $X$ be a normal toric variety such that $\Sigma$ has full-dimensional convex support. We define the \emph{moving cone} of the secondary fan of $X$ to be
    \begin{align*}
        \Mov_{\text{GKZ}} \equ \bigcup_{\Sigma: I_\emptyset = \emptyset} \Gamma_{\Sigma, \emptyset},
    \end{align*}
    i.e.\ the union of all GKZ cones with $I_\emptyset = \emptyset$.
\end{definition}

For toric varieties, the above description coincides with the \emph{moving cone} of $X$, which is a classical object in birational geometry.

\begin{definition}\label{def:movable_div}
    Let $X_\Sigma$ be a normal toric variety such that $\Sigma$ has full-dimensional convex support.
    An effective divisor $D \ge 0$ on $X$ has a \emph{fixed component}, if there is an effective divisor $D_0 \neq 0$ such that every effective divisor $D'\sim D$ satisfies $D' \ge D_0$. Hence, a Weil divisor $D$ on $X$ is called \emph{movable} if it does not have a fixed component, and the \emph{moving cone of $X$} is defined by 
    \begin{align*}
        \Mov(X) \deq \overline{\Cone}([D] \mid D \text{ is Cartier and movable}\} \subseteq N^1(X).
    \end{align*}
\end{definition}

Now as $\Proj^D(S)$ is always simplicial in this setting (cf.\ \cite{BS}, Proposition 3.4), we can apply \cite{CLS}, Theorem 15.1.10, and deduce that both $\Mov_{\text{GKZ}}$ and $\Mov(X)$ coincide. In particular, for $X = X_\Sigma$ it holds
\begin{align*}
    \Mov(X_\Sigma) \equ \bigcup_{\Sigma'} \Nef(X_{\Sigma'}),
\end{align*}
where the union is over all simplicial fans $\Sigma'$ in $N_\BR$ such that $X_{\Sigma'}$ is semiprojective and $\Sigma'(1) = \Sigma(1)$.

The moving cone can be described by relevant elements in the following way:

\begin{proposition}\label{prop:movable_cone_relevant}
    Let $X_\Sigma$ be a quasiprojective simplicial normal toric variety such that $\Sigma$ has full-dimensional convex support. Then $S = \Cox(X_\Sigma)$ is a polynomial ring, say $S = k[x_1, \ldots, x_n]$. We assume that $r = \rank(\Cl(X)) < n$. Then the movable cone of the secondary fan is given by
    \begin{align*}
        \Mov_{\text{GKZ}} \equ \bigcap_{i=1}^n \CC_D(f_i),
    \end{align*}
    where the $f_i$ are given by a product of $n-1$ variables $f_i = x_1 \cdot \ldots x_{i-1} \cdot x_{i+1} \cdot  \ldots \cdot x_n$ for $i = 1, \ldots n$, i.e.\ $f_i$ omits $x_i$.
\end{proposition}

\begin{proof}
    The statement essentially follows from \cite{CLS}, Proposition 15.2.4, stating that
        \begin{align*}
        \Mov_{\text{GKZ}} \equ \bigcap_{i=1}^n \overline{\Cone}(\beta_1, \ldots, \widehat{\beta_i}, \ldots, \beta_n),
    \end{align*}
    where $\widehat{\beta_i}$ means that $\beta_i$ is ommited. The reason for that is that for some $\chi \otimes 1 \in C_\beta$ it holds
    \begin{align*}
        F_{i,\chi} \neq \emptyset \ \iff \ \chi \otimes 1 \in \overline{\Cone}(\beta_1, \ldots, \widehat{\beta_i}, \ldots, \beta_n)
    \end{align*}
    and
    \begin{align*}
        I_\emptyset = \emptyset \ \iff \ \chi \otimes 1 \in \bigcap_{i=1}^n\overline{\Cone}(\beta_1, \ldots, \widehat{\beta_i}, \ldots, \beta_n).
    \end{align*}    
    Now clearly $\deg(x_i) = \beta_i$ by construction and the claim follows by Definition~\ref{def:moving_GKZ_cone}.
\end{proof}

\begin{example}\label{ref:standard_ex_moving}
Consider the toric variety $X = \Bl_p(\BP^2)$ from the previous examples. By Proposition~\ref{prop:movable_cone_relevant}, we have to take all products of exactly three variables into account, so that
    \begin{align*}
        \Mov_{\text{GKZ}} \equ \CC_D(xyz) \cap \CC_D(xyw) \cap \CC_D(yzw) \cap \CC_D(xzw) \equ \CC_D(xz)
    \end{align*}
    exactly gives the chamber $\overline{\Cone}((1,0), (1,1))$, which gives the section ring $S' = \BC[x,y,z, xw, yw]$. We compute $\Cl(X)$: Consider $\BP^2$ with homogeneous coordinates $[x:y:z]$. Then $X$ contains the (class of) the strict transform of the line $(z=0)$, say $H$, and the exceptional divisor $E$, say corresponding to $w$. Then both $x$ and $y$ correspond to the class of $H-E$, as the divisors $D_x$ and $D_y$ associated to $x$ and $y$ are linearly equivalent. Then, as $E$ has a negative self-intersection number, it holds that 
    \begin{align*}
        \Nef(X) \equ \overline{\Cone}([H], [H-E]) \equ \overline{\Cone}((1,0), (1,1)) \equ \Mov_{\text{GKZ}}.
    \end{align*}
\end{example}


\subsection{Secondary Fans for Toric Prevarieties}

In Example~\ref{ex:standard_ex_chambers}, we were able to compute a secondary fan for a nonseparated toric variety.
Actually, none of the construction steps in computing the secondary fan for $X_\Sigma$, being a quasiprojective simplicial normal toric variety such that $\Sigma$ has full-dimensional convex support, relies on $X$ being separated. In fact, the computations depend on the Cox ring of $X$ and its grading only, and it is not hard to see that Cox rings of toric prevarieties are just polynomial ring with finitely many variables, as is the case for toric varieties.
In particular, the class group of a toric prevariety remains to be generated by the classes of the torus invariant prime divisors (cf.\ \cite{CLS}, Theorem 4.0.20).
Note that the key reason for this is that we can describe the chambers in the secondary fan purely in terms of subsets of the irrelevant ideal $\Cox(X)_+$ (in the sense of \cite{BS}). Therefore, we can state:

\begin{theorem}\label{thm:secondary_prev}
    Let $X_\Sigma$ be a quasiprojective simplicial normal toric prevariety such that $\Sigma$ has full-dimensional convex support, and $S = \Cox(X)$. Then $\Cl(X)_\BR$ has a decomposition into chambers 
    \begin{align*}
           \sigma_a \equ \bigcap_{\substack{f \in \Gen^D(S) \\ a \in \CC_D(f)}} \CC_D(f).
    \end{align*}
    for $a \in \lint(\sigma(S))$. In particular, Theorem~\ref{thm:list_equiv_chamber_var} remains valid for simplicial toric prevarieties.
\end{theorem}

\begin{proof}
    \cite{paper1}, Lemma 1.19 implies that the definition of $B(\chi)$ does not change when transitioning to a nonseparated toric variety, as this purely depends on the structure of the Cox ring. Thus, we can transport the notion of generic characters and their properties to the nonseparated world. The claim follows.
\end{proof}


While the construction of GKZ cones seems not to be immediately generalizable, the description in Proposition~\ref{prop:movable_cone_relevant} is. In fact, it is taken as a definition in \cite{CRB}.  

\begin{definition}\label{def:moving_cone_CRB}
    Let $\CW = (w_1, \ldots, w_n)$ be a family of vectors in $D_\BQ$ that generate $D_\BQ$. The \emph{moving cone} of $\CW$ is 
    \begin{align*}
        \Mov(\CW) \deq \bigcap_{i=1}^n \overline{\Cone}(w_j \mid j \neq i) \subseteq D_\BQ.
    \end{align*}
\end{definition}

One can see that the above Definition generalizes Definition~\ref{def:moving_GKZ_cone}, so that we may speak about the secondary fan and the moving cone of a toric prevariety. For a list of properties of $\Mov(\CW)$, we refer to \cite{CRB}, Theorem 2.2.2.6.

\begin{example}\label{ex:standard_ex_mov_prev}
Let again $S = \BC[x, y, z, w]$ and $D = \BZ^2$, where $\deg(x) = \deg(y) = (1, 0)$, $\deg(z) = (1, 1)$ and $\deg(w) = (0, 1)$ (cf.\ Example~\ref{ex:double_origin}). This time we take $X = \Proj^D(S)$. Then by Definition~\ref{def:moving_cone_CRB} the moving cone of the vector configuration $(e_1, e_1, e_1+e_2, e_2)$ is given by
    \begin{align*}
        \Mov_{\text{GKZ}} &\equ \overline{\Cone}(e_1, e_1+e_2) \cap \overline{\Cone}(e_1, e_2) \cap \overline{\Cone}(e_1, e_1+e_2, e_2)  \\ &\equ \overline{\Cone}(e_1, e_1+e_2) \\ &\equ \CC_D(xz),
    \end{align*}
    as expected. In particular, the secondary fan and the moving cone of the toric prevariety $X$ are the same as for $\Bl_p(\BP^2)$. The reason for this is, of course, that $X$ is the gluing of the VGIT quotients of $\Bl_p(\BP^2)$.
\end{example}

\section{Birational Geometry via $D$-graded Proj}\label{sec:MDS}
In this short section, we generalize the insights from §3 on simplicial normal toric (pre)varieties to Mori dream spaces. In particular, we deduce an explicit form of Mori equivalence that can be generalized to more general varieties, such as Calabi-Yau varieties.

\begin{theorem}\label{thm:MDS_relevant} 
    Let $X$ be a Mori dream space, $S = \Cox(X)$ and $C \in \Cl(X)$ a ($\BQ$-Cartier) divisor satisfying $\deg(C) \in \lint(\sigma(S))$. Then the Mori chamber containing $\deg(C)$ is given by
    \begin{align*}
           \sigma_C \equ \bigcap_{\substack{f \in \Gen^D(S) \\ \deg(C) \in \CC_D(f)}} \CC_D(f).
    \end{align*}
\end{theorem}

\begin{proof}
    By \cite{HK}, Proposition 2.9, the Mori chambers correspond exactly to the GIT-cones, that is, the chambers of the weight space $\Eff(X)_\BR$. According to \cite{DH}, Theorem 3.3.2, the interior of a maximal-dimensional GIT-cone is exactly the area where $X^\text{ss}(d) = X^{s}(d)$ for all $d$ in that interior. Hence, the notion of the GIT-cone coincides with \cite{CRB} and the statement follows from Corollary~\ref{cor:GIT_relevant_cones}.

\end{proof}

As in §3, each Mori chamber encodes a GIT quotient associated to a linearized line bundle on $X$. The explicit formula for the Mori chamber from Theorem~\ref{thm:MDS_relevant} directly translates to Mori equivalence.

\begin{corollary}\label{cor:MDS_chamber=sectionring}
    Let $X$ be a Mori dream space, $S = \Cox(X)$ and $C_1, C_2 \in \Cl(X)$. Then $C_1$ and $C_2$ are Mori equivalent if and only if for all relevant $f \in S$ it holds that
    \begin{align*}
        \deg(C_1) \in \CC_D(f) \ \ \iff \ \ \deg(C_2) \in \CC_D(f).
    \end{align*}
    In particular, the variety associated with the chamber $C$ is precisely the multigraded Proj of the section ring supported in that chamber.
\end{corollary}

\begin{remark}
    We expect that this chamber structure holds true for any multigraded ring $S$. Concretely, the VGIT of the action of $G$ on $\Spec(S)$ correspond to the chamber decomposition of $\sigma(S)$, where the chamber containing $a \in \lint(\sigma(S))$ is given by
    \begin{align*}
                   \sigma_a \equ \bigcap_{\substack{f \in \Rel^D(S) \\ a \in \CC_D(f)}} \CC_D(f).
    \end{align*}
\end{remark}

\end{document}